\documentclass[a4paper,leqno]{article}
\usepackage{amsfonts}
\usepackage{amsmath}
\usepackage{amssymb}
\numberwithin{equation}{subsection}

\makeatletter
\def\diagram{\leftwidth=\z@ \rightwidth=\z@ \topheight=\z@
\botheight=\z@ \setbox\@picbox\hbox\bgroup}
\def\enddiagram{\egroup\wd\@picbox\rightwidth\unitlength
\ht\@picbox\topheight\unitlength \dp\@picbox\botheight\unitlength
\hskip\leftwidth\unitlength\box\@picbox}
\def\bfig{\begin{diagram}}
\def\efig{\end{diagram}}
\newcount\wideness \newcount\leftwidth \newcount\rightwidth
\newcount\highness \newcount\topheight \newcount\botheight
\def\ratchet#1#2{\ifnum#1<#2 \global #1=#2 \fi}
\def\putbox(#1,#2)#3{%
\horsize{\wideness}{#3} \divide\wideness by 2 {\advance\wideness
by #1 \ratchet{\rightwidth}{\wideness}} {\advance\wideness by -#1
\ratchet{\leftwidth}{\wideness}} \vertsize{\highness}{#3}
\divide\highness by 2 {\advance\highness by #2
\ratchet{\topheight}{\highness}} {\advance\highness by -#2
\ratchet{\botheight}{\highness}} \put(#1,#2){\makebox(0,0){$#3$}}}
\def\putlbox(#1,#2)#3{%
\horsize{\wideness}{#3} {\advance\wideness by #1
\ratchet{\rightwidth}{\wideness}} {\ratchet{\leftwidth}{-#1}}
\vertsize{\highness}{#3} \divide\highness by 2 {\advance\highness
by #2 \ratchet{\topheight}{\highness}} {\advance\highness by -#2
\ratchet{\botheight}{\highness}}
\put(#1,#2){\makebox(0,0)[l]{$#3$}}}
\def\putrbox(#1,#2)#3{%
\horsize{\wideness}{#3} {\ratchet{\rightwidth}{#1}}
{\advance\wideness by -#1 \ratchet{\leftwidth}{\wideness}}
\vertsize{\highness}{#3} \divide\highness by 2 {\advance\highness
by #2 \ratchet{\topheight}{\highness}} {\advance\highness by -#2
\ratchet{\botheight}{\highness}}
\put(#1,#2){\makebox(0,0)[r]{$#3$}}}

\def\adjust[#1]{} 
\newcount \coefa
\newcount \coefb
\newcount \coefc
\newcount\tempcounta
\newcount\tempcountb
\newcount\tempcountc
\newcount\tempcountd
\newcount\xext
\newcount\yext
\newcount\xoff
\newcount\yoff
\newcount\gap%
\newcount\arrowtypea
\newcount\arrowtypeb
\newcount\arrowtypec
\newcount\arrowtyped
\newcount\arrowtypee
\newcount\height
\newcount\width
\newcount\xpos
\newcount\ypos
\newcount\run
\newcount\rise
\newcount\arrowlength
\newcount\halflength
\newcount\arrowtype
\newdimen\tempdimen
\newdimen\xlen
\newdimen\ylen
\newsavebox{\tempboxa}%
\newsavebox{\tempboxb}%
\newsavebox{\tempboxc}%
\newdimen\w@dth
\def\setw@dth#1#2{\setbox\z@\hbox{$#1$}\w@dth=\wd\z@
\setbox\@ne\hbox{$#2$}\ifnum\w@dth<\wd\@ne \w@dth=\wd\@ne \fi
\advance\w@dth by 1.2em}
\def\t@^#1_#2{\def\n@one{#1}\def\n@two{#2}\mathrel{\setw@dth{#1}{#2}
\mathop{\hbox to \w@dth{\rightarrowfill}}\limits
\ifx\n@one\empty\else ^{\box\z@}\fi \ifx\n@two\empty\else
_{\box\@ne}\fi}}
\def\t@@^#1{\@ifnextchar_ {\t@^{#1}}{\t@^{#1}_{}}}
\def\to{\@ifnextchar^ {\t@@}{\t@@^{}}}
\def\t@left^#1_#2{\def\n@one{#1}\def\n@two{#2}\mathrel{\setw@dth{#1}{#2}
\mathop{\hbox to \w@dth{\leftarrowfill}}\limits
\ifx\n@one\empty\else ^{\box\z@}\fi \ifx\n@two\empty\else
_{\box\@ne}\fi}}
\def\t@@left^#1{\@ifnextchar_ {\t@left^{#1}}{\t@left^{#1}_{}}}
\def\toleft{\@ifnextchar^ {\t@@left}{\t@@left^{}}}
\def\two@^#1_#2{\def\n@one{#1}\def\n@two{#2}\mathrel{\setw@dth{#1}{#2}
\mathop{\vcenter{\hbox to \w@dth{\rightarrowfill}\kern-1.7ex
                 \hbox to \w@dth{\rightarrowfill}}%
       }\limits
\ifx\n@one\empty\else ^{\box\z@}\fi \ifx\n@two\empty\else
_{\box\@ne}\fi}}
\def\tw@@^#1{\@ifnextchar_ {\two@^{#1}}{\two@^{#1}_{}}}
\def\two{\@ifnextchar^ {\tw@@}{\tw@@^{}}}
\def\tofr@^#1_#2{\def\n@one{#1}\def\n@two{#2}\mathrel{\setw@dth{#1}{#2}
\mathop{\vcenter{\hbox to \w@dth{\rightarrowfill}\kern-1.7ex
                 \hbox to \w@dth{\leftarrowfill}}%
       }\limits
\ifx\n@one\empty\else ^{\box\z@}\fi \ifx\n@two\empty\else
_{\box\@ne}\fi}}
\def\t@fr@^#1{\@ifnextchar_ {\tofr@^{#1}}{\tofr@^{#1}_{}}}
\def\tofro{\@ifnextchar^ {\t@fr@}{\t@fr@^{}}}

\def\mon{\mathop{\m@th\hbox to
      14.6\P@{\lasyb\char'51\hskip-2.1\P@$\arrext$\hss
$\mathord\rightarrow$}}\limits} 
\def\leftmono{\mathrel{\m@th\hbox to
14.6\P@{$\mathord\leftarrow$\hss$\arrext$\hskip-2.1\P@\lasyb\char'50%
}}\limits} 
\mathchardef\arrext="0200
\def\settypes(#1,#2,#3){\arrowtypea#1 \arrowtypeb#2 \arrowtypec#3}
\def\settoheight#1#2{\setbox\@tempboxa\hbox{#2}#1\ht\@tempboxa\relax}%
\def\settodepth#1#2{\setbox\@tempboxa\hbox{#2}#1\dp\@tempboxa\relax}%
\def\settokens[#1`#2`#3`#4]{%
     \def\tokena{#1}\def\tokenb{#2}\def\tokenc{#3}\def\tokend{#4}}
\def\setsqparms[#1`#2`#3`#4;#5`#6]{%
\arrowtypea #1 \arrowtypeb #2 \arrowtypec #3 \arrowtyped #4
\width #5 \height #6 }
\def\setpos(#1,#2){\xpos=#1 \ypos#2}

\def\settriparms[#1`#2`#3;#4]{\settripairparms[#1`#2`#3`1`1;#4]}%
\def\settripairparms[#1`#2`#3`#4`#5;#6]{%
\arrowtypea #1 \arrowtypeb #2 \arrowtypec #3 \arrowtyped #4
\arrowtypee #5 \width #6 \height #6 }
\def\resetparms{\settripairparms[1`1`1`1`1;500]\width 500}
\resetparms
\def\mvector(#1,#2)#3{
\put(0,0){\vector(#1,#2){#3}}%
\put(0,0){\vector(#1,#2){26}}%
}
\def\evector(#1,#2)#3{{
\arrowlength #3
\put(0,0){\vector(#1,#2){\arrowlength}}%
\advance \arrowlength by-30
\put(0,0){\vector(#1,#2){\arrowlength}}%
}}
\def\horsize#1#2{%
\settowidth{\tempdimen}{$#2$}%
#1=\tempdimen \divide #1 by\unitlength }
\def\vertsize#1#2{%
\settoheight{\tempdimen}{$#2$}%
#1=\tempdimen
\settodepth{\tempdimen}{$#2$}%
\advance #1 by\tempdimen \divide #1 by\unitlength }
\def\putvector(#1,#2)(#3,#4)#5#6{{%
\ifnum3<\arrowtype \putdashvector(#1,#2)(#3,#4)#5\arrowtype \else
\ifnum\arrowtype<-3 \putdashvector(#1,#2)(#3,#4)#5\arrowtype \else
\xpos=#1 \ypos=#2 \run=#3 \rise=#4 \arrowlength=#5 \ifnum
\arrowtype<0
    \ifnum \run=0
        \advance \ypos by-\arrowlength
    \else
        \tempcounta \arrowlength
        \multiply \tempcounta by\rise
        \divide \tempcounta by\run
        \ifnum\run>0
            \advance \xpos by\arrowlength
            \advance \ypos by\tempcounta
        \else
            \advance \xpos by-\arrowlength
            \advance \ypos by-\tempcounta
        \fi
    \fi
    \multiply \arrowtype by-1
    \multiply \rise by-1
    \multiply \run by-1
\fi \ifcase \arrowtype
\or \put(\xpos,\ypos){\vector(\run,\rise){\arrowlength}}%
\or \put(\xpos,\ypos){\mvector(\run,\rise)\arrowlength}%
\or \put(\xpos,\ypos){\evector(\run,\rise){\arrowlength}}%
\fi\fi\fi }}
\def\putsplitvector(#1,#2)#3#4{
\xpos #1 \ypos #2 \arrowtype #4 \halflength #3 \arrowlength #3
\gap 140 \advance \halflength by-\gap \divide \halflength by2
\ifnum\arrowtype>0
   \ifcase \arrowtype
   \or \put(\xpos,\ypos){\line(0,-1){\halflength}}%
       \advance\ypos by-\halflength
       \advance\ypos by-\gap
       \put(\xpos,\ypos){\vector(0,-1){\halflength}}%
   \or \put(\xpos,\ypos){\line(0,-1)\halflength}%
       \put(\xpos,\ypos){\vector(0,-1)3}%
       \advance\ypos by-\halflength
       \advance\ypos by-\gap
       \put(\xpos,\ypos){\vector(0,-1){\halflength}}%
   \or \put(\xpos,\ypos){\line(0,-1)\halflength}%
       \advance\ypos by-\halflength
       \advance\ypos by-\gap
       \put(\xpos,\ypos){\evector(0,-1){\halflength}}%
   \fi
\else \arrowtype=-\arrowtype
   \ifcase\arrowtype
   \or \advance \ypos by-\arrowlength
       \put(\xpos,\ypos){\line(0,1){\halflength}}%
       \advance\ypos by\halflength
       \advance\ypos by\gap
       \put(\xpos,\ypos){\vector(0,1){\halflength}}%
   \or \advance \ypos by-\arrowlength
       \put(\xpos,\ypos){\line(0,1)\halflength}%
       \put(\xpos,\ypos){\vector(0,1)3}%
       \advance\ypos by\halflength
       \advance\ypos by\gap
       \put(\xpos,\ypos){\vector(0,1){\halflength}}%
   \or \advance \ypos by-\arrowlength
       \put(\xpos,\ypos){\line(0,1)\halflength}%
       \advance\ypos by\halflength
       \advance\ypos by\gap
       \put(\xpos,\ypos){\evector(0,1){\halflength}}%
   \fi
\fi }
\def\putmorphism(#1)(#2,#3)[#4`#5`#6]#7#8#9{{%
\run #2 \rise #3 \ifnum\rise=0
  \puthmorphism(#1)[#4`#5`#6]{#7}{#8}#9%
\else\ifnum\run=0
  \putvmorphism(#1)[#4`#5`#6]{#7}{#8}#9%
\else
\setpos(#1)%
\arrowlength #7 \arrowtype #8 \ifnum\run=0 \else\ifnum\rise=0
\else \ifnum\run>0
    \coefa=1
\else
   \coefa=-1
\fi \ifnum\arrowtype>0
   \coefb=0
   \coefc=-1
\else
   \coefb=\coefa
   \coefc=1
   \arrowtype=-\arrowtype
\fi \width=2 \multiply \width by\run \divide \width by\rise
\ifnum \width<0  \width=-\width\fi \advance\width by60 \if l#9
\width=-\width\fi
\putbox(\xpos,\ypos){#4}
{\multiply \coefa by\arrowlength
\advance\xpos by\coefa \multiply \coefa by\rise \divide \coefa
by\run \advance \ypos by\coefa
\putbox(\xpos,\ypos){#5} }%
{\multiply \coefa by\arrowlength
\divide \coefa by2 \advance \xpos by\coefa \advance \xpos by\width
\multiply \coefa by\rise \divide \coefa by\run \advance \ypos
by\coefa
\if l#9%
   \putrbox(\xpos,\ypos){#6}%
\else\if r#9%
   \putlbox(\xpos,\ypos){#6}%
\fi\fi }%
{\multiply \rise by-\coefc
\multiply \run by-\coefc \multiply \coefb by\arrowlength \advance
\xpos by\coefb \multiply \coefb by\rise \divide \coefb by\run
\advance \ypos by\coefb \multiply \coefc by70 \advance \ypos
by\coefc \multiply \coefc by\run \divide \coefc by\rise \advance
\xpos by\coefc \multiply \coefa by140 \multiply \coefa by\run
\divide \coefa by\rise \advance \arrowlength by\coefa
\ifcase\arrowtype
\or \put(\xpos,\ypos){\vector(\run,\rise){\arrowlength}}%
\or \put(\xpos,\ypos){\mvector(\run,\rise){\arrowlength}}%
\or \put(\xpos,\ypos){\evector(\run,\rise){\arrowlength}}%
\fi}\fi\fi\fi\fi}}

\newcount\numbdashes \newcount\lengthdash \newcount\increment
\def\howmanydashes{
\numbdashes=\arrowlength \lengthdash=40 \divide\numbdashes by
\lengthdash \lengthdash=\arrowlength \divide\lengthdash by
\numbdashes
\increment=\lengthdash \multiply\lengthdash by 3
\divide\lengthdash by 5 }
\def\putdashvector(#1)(#2,#3)#4#5{%
\ifnum#3=0 \putdashhvector(#1){#4}#5 \else \ifnum#2=0
\putdashvvector(#1){#4}#5\fi\fi}
\def\putdashhvector(#1,#2)#3#4{{%
\arrowlength=#3 \howmanydashes
\multiput(#1,#2)(\increment,0){\numbdashes}%
{\vrule height .4pt width \lengthdash\unitlength} \arrowtype=#4
\xpos=#1 \ifnum\arrowtype<0 \advance\arrowtype by 7 \fi
\ifcase\arrowtype \or \advance\xpos by 10
    \put(\xpos,#2){\vector(-1,0){\lengthdash}}
    \advance\xpos by 40
    \put(\xpos,#2){\vector(-1,0){\lengthdash}}
\or \advance \xpos by 10
    \put(\xpos,#2){\vector(-1,0){\lengthdash}}
    \advance\xpos by  \arrowlength
    \advance\xpos by  -50
    \put(\xpos,#2){\vector(-1,0){\lengthdash}}
\or \advance\xpos by 10
    \put(\xpos,#2){\vector(-1,0){\lengthdash}}
\or \advance\xpos by \arrowlength
    \advance\xpos by -\lengthdash
    \put(\xpos,#2){\vector(1,0){\lengthdash}}
\or {\advance\xpos by 10
    \put(\xpos,#2){\vector(1,0){\lengthdash}}}
    \advance\xpos by \arrowlength
    \advance\xpos by -\lengthdash
    \put(\xpos,#2){\vector(1,0){\lengthdash}}
\or \advance\xpos by \arrowlength
    \advance\xpos by -\lengthdash
    \put(\xpos,#2){\vector(1,0){\lengthdash}}
    \advance\xpos by -40
    \put(\xpos,#2){\vector(1,0){\lengthdash}}
   \fi
}}
\def\putdashvvector(#1,#2)#3#4{{%
\arrowlength=#3 \howmanydashes \ypos=#2 \advance\ypos by
-\arrowlength
\multiput(#1,#2)(0,\increment){\numbdashes}%
    {\vrule width .4pt height \lengthdash\unitlength}
\arrowtype=#4 \ypos=#2 \ifnum\arrowtype<0 \advance\arrowtype by 7
\fi \ifcase\arrowtype \or \advance\ypos by \arrowlength
\advance\ypos by -40
    \put(#1,\ypos){\vector(0,1){\lengthdash}}
    \advance\ypos by -40
    \put(#1,\ypos){\vector(0,1){\lengthdash}}
\or \advance\ypos by 10
    \put(#1,\ypos){\vector(0,1){\lengthdash}}
    \advance\ypos by \arrowlength \advance\ypos by -40
    \put(#1,\ypos){\vector(0,1){\lengthdash}}
\or \advance\ypos by \arrowlength \advance\ypos by -40
    \put(#1,\ypos){\vector(0,1){\lengthdash}}
\or \advance\ypos by 10
    \put(#1,\ypos){\vector(0,-1){\lengthdash}}
\or \advance\ypos by 10
    \put(#1,\ypos){\vector(0,-1){\lengthdash}}
    \advance\ypos by \arrowlength \advance\ypos by -40
    \put(#1,\ypos){\vector(0,-1){\lengthdash}}
\or \advance\ypos by 10
    \put(#1,\ypos){\vector(0,-1){\lengthdash}}
    \advance\ypos by 40
    \put(#1,\ypos){\vector(0,-1){\lengthdash}}
\fi }}
\def\puthmorphism(#1,#2)[#3`#4`#5]#6#7#8{{%
\xpos #1 \ypos #2 \width #6 \arrowlength #6 \arrowtype=#7
\putbox(\xpos,\ypos){#3\vphantom{#4}}%
{\advance \xpos by\arrowlength
\putbox(\xpos,\ypos){\vphantom{#3}#4}}%
\horsize{\tempcounta}{#3}%
\horsize{\tempcountb}{#4}%
\divide \tempcounta by2 \divide \tempcountb by2 \advance
\tempcounta by30 \advance \tempcountb by30 \advance \xpos
by\tempcounta \advance \arrowlength by-\tempcounta \advance
\arrowlength by-\tempcountb
\putvector(\xpos,\ypos)(1,0)\arrowlength\arrowtype \divide
\arrowlength by2 \advance \xpos by\arrowlength
\vertsize{\tempcounta}{#5}%
\divide\tempcounta by2 \advance \tempcounta by20
\if a#8 %
   \advance \ypos by\tempcounta
   \putbox(\xpos,\ypos){#5}%
\else
   \advance \ypos by-\tempcounta
   \putbox(\xpos,\ypos){#5}%
\fi}}
\def\putvmorphism(#1,#2)[#3`#4`#5]#6#7#8{{%
\xpos #1 \ypos #2 \arrowlength #6 \arrowtype #7
\settowidth{\xlen}{$#5$}%
\putbox(\xpos,\ypos){#3}%
{\advance \ypos by-\arrowlength
\putbox(\xpos,\ypos){#4}}%
{\advance\arrowlength by-140 \advance \ypos by-70 \ifdim\xlen>0pt
   \if m#8%
      \putsplitvector(\xpos,\ypos)\arrowlength\arrowtype
   \else
   \putvector(\xpos,\ypos)(0,-1)\arrowlength\arrowtype
   \fi
\else
   \putvector(\xpos,\ypos)(0,-1)\arrowlength\arrowtype
\fi}%
\ifdim\xlen>0pt
   \divide \arrowlength by2
   \advance\ypos by-\arrowlength
   \if l#8%
      \advance \xpos by-40
      \putrbox(\xpos,\ypos){#5}%
   \else\if r#8%
      \advance \xpos by40
      \putlbox(\xpos,\ypos){#5}%
   \else
      \putbox(\xpos,\ypos){#5}%
   \fi\fi
\fi }}
\def\putsquarep<#1>(#2)[#3;#4`#5`#6`#7]{{%
\setsqparms[#1]%
\setpos(#2)%
\settokens[#3]%
\puthmorphism(\xpos,\ypos)[\tokenc`\tokend`{#7}]{\width}{\arrowtyped}b%
\advance\ypos by \height
\puthmorphism(\xpos,\ypos)[\tokena`\tokenb`{#4}]{\width}{\arrowtypea}a%
\putvmorphism(\xpos,\ypos)[``{#5}]{\height}{\arrowtypeb}l%
\advance\xpos by \width
\putvmorphism(\xpos,\ypos)[``{#6}]{\height}{\arrowtypec}r%
}}
\def\putsquare{\@ifnextchar <{\putsquarep}{\putsquarep%
   <\arrowtypea`\arrowtypeb`\arrowtypec`\arrowtyped;\width`\height>}}
\def\square{\@ifnextchar< {\squarep}{\squarep
   <\arrowtypea`\arrowtypeb`\arrowtypec`\arrowtyped;\width`\height>}}
\def\squarep<#1>[#2`#3`#4`#5;#6`#7`#8`#9]{{
\setsqparms[#1]
\diagram
\putsquarep<\arrowtypea`\arrowtypeb`\arrowtypec`
\arrowtyped;\width`\height>
(0,0)[#2`#3`#4`{#5};#6`#7`#8`{#9}]
\enddiagram
}}                                                 
\def\putptrianglep<#1>(#2,#3)[#4`#5`#6;#7`#8`#9]{{%
\settriparms[#1]%
\xpos=#2 \ypos=#3 \advance\ypos by \height
\puthmorphism(\xpos,\ypos)[#4`#5`{#7}]{\height}{\arrowtypea}a%
\putvmorphism(\xpos,\ypos)[`#6`{#8}]{\height}{\arrowtypeb}l%
\advance\xpos by\height
\putmorphism(\xpos,\ypos)(-1,-1)[``{#9}]{\height}{\arrowtypec}r%
}}
\def\putptriangle{\@ifnextchar <{\putptrianglep}{\putptrianglep
   <\arrowtypea`\arrowtypeb`\arrowtypec;\height>}}
\def\ptriangle{\@ifnextchar <{\ptrianglep}{\ptrianglep
   <\arrowtypea`\arrowtypeb`\arrowtypec;\height>}}
\def\ptrianglep<#1>[#2`#3`#4;#5`#6`#7]{{
\settriparms[#1]
\diagram
\putptrianglep<\arrowtypea`\arrowtypeb`
\arrowtypec;\height>
(0,0)[#2`#3`#4;#5`#6`{#7}]
\enddiagram
}}                                            
\def\putqtrianglep<#1>(#2,#3)[#4`#5`#6;#7`#8`#9]{{%
\settriparms[#1]%
\xpos=#2 \ypos=#3 \advance\ypos by\height
\puthmorphism(\xpos,\ypos)[#4`#5`{#7}]{\height}{\arrowtypea}a%
\putmorphism(\xpos,\ypos)(1,-1)[``{#8}]{\height}{\arrowtypeb}l%
\advance\xpos by\height
\putvmorphism(\xpos,\ypos)[`#6`{#9}]{\height}{\arrowtypec}r%
}}
\def\putqtriangle{\@ifnextchar <{\putqtrianglep}{\putqtrianglep
   <\arrowtypea`\arrowtypeb`\arrowtypec;\height>}}
\def\qtriangle{\@ifnextchar <{\qtrianglep}{\qtrianglep
   <\arrowtypea`\arrowtypeb`\arrowtypec;\height>}}
\def\qtrianglep<#1>[#2`#3`#4;#5`#6`#7]{{
\settriparms[#1]
\width=\height                                
\diagram
\putqtrianglep<\arrowtypea`\arrowtypeb`
\arrowtypec;\height>
(0,0)[#2`#3`#4;#5`#6`{#7}]
\enddiagram
}}
\def\putdtrianglep<#1>(#2,#3)[#4`#5`#6;#7`#8`#9]{{%
\settriparms[#1]%
\xpos=#2 \ypos=#3
\puthmorphism(\xpos,\ypos)[#5`#6`{#9}]{\height}{\arrowtypec}b%
\advance\xpos by \height \advance\ypos by\height
\putmorphism(\xpos,\ypos)(-1,-1)[``{#7}]{\height}{\arrowtypea}l%
\putvmorphism(\xpos,\ypos)[#4``{#8}]{\height}{\arrowtypeb}r%
}}
\def\putdtriangle{\@ifnextchar <{\putdtrianglep}{\putdtrianglep
   <\arrowtypea`\arrowtypeb`\arrowtypec;\height>}}
\def\dtriangle{\@ifnextchar <{\dtrianglep}{\dtrianglep
   <\arrowtypea`\arrowtypeb`\arrowtypec;\height>}}
\def\dtrianglep<#1>[#2`#3`#4;#5`#6`#7]{{
\settriparms[#1]
\width=\height                                
\diagram
\putdtrianglep<\arrowtypea`\arrowtypeb`
\arrowtypec;\height>
(0,0)[#2`#3`#4;#5`#6`{#7}]
\enddiagram
}}
\def\putbtrianglep<#1>(#2,#3)[#4`#5`#6;#7`#8`#9]{{%
\settriparms[#1]%
\xpos=#2 \ypos=#3
\puthmorphism(\xpos,\ypos)[#5`#6`{#9}]{\height}{\arrowtypec}b%
\advance\ypos by\height
\putmorphism(\xpos,\ypos)(1,-1)[``{#8}]{\height}{\arrowtypeb}r%
\putvmorphism(\xpos,\ypos)[#4``{#7}]{\height}{\arrowtypea}l%
}}
\def\putbtriangle{\@ifnextchar <{\putbtrianglep}{\putbtrianglep
   <\arrowtypea`\arrowtypeb`\arrowtypec;\height>}}
\def\btriangle{\@ifnextchar <{\btrianglep}{\btrianglep
   <\arrowtypea`\arrowtypeb`\arrowtypec;\height>}}
\def\btrianglep<#1>[#2`#3`#4;#5`#6`#7]{{
\settriparms[#1]
\width=\height                               
\diagram
\putbtrianglep<\arrowtypea`\arrowtypeb`
\arrowtypec;\height>
(0,0)[#2`#3`#4;#5`#6`{#7}]
\enddiagram
}}
\def\putAtrianglep<#1>(#2,#3)[#4`#5`#6;#7`#8`#9]{{%
\settriparms[#1]%
\xpos=#2 \ypos=#3 {\multiply \height by2
\puthmorphism(\xpos,\ypos)[#5`#6`{#9}]{\height}{\arrowtypec}b}%
\advance\xpos by\height \advance\ypos by\height
\putmorphism(\xpos,\ypos)(-1,-1)[#4``{#7}]{\height}{\arrowtypea}l%
\putmorphism(\xpos,\ypos)(1,-1)[``{#8}]{\height}{\arrowtypeb}r%
}}
\def\putAtriangle{\@ifnextchar <{\putAtrianglep}{\putAtrianglep
   <\arrowtypea`\arrowtypeb`\arrowtypec;\height>}}
\def\Atriangle{\@ifnextchar <{\Atrianglep}{\Atrianglep
   <\arrowtypea`\arrowtypeb`\arrowtypec;\height>}}
\def\Atrianglep<#1>[#2`#3`#4;#5`#6`#7]{{
\settriparms[#1]
\width=\height                                     
\diagram
\putAtrianglep<\arrowtypea`\arrowtypeb`
\arrowtypec;\height>
(0,0)[#2`#3`#4;#5`#6`{#7}]
\enddiagram
}}
\def\putAtrianglepairp<#1>(#2)[#3;#4`#5`#6`#7`#8]{{%
\settripairparms[#1]%
\setpos(#2)%
\settokens[#3]%
\puthmorphism(\xpos,\ypos)[\tokenb`\tokenc`{#7}]{\height}{\arrowtyped}b%
\advance\xpos by\height
\puthmorphism(\xpos,\ypos)[\phantom{\tokenc}`\tokend`{#8}]%
{\height}{\arrowtypee}b%
\advance\ypos by\height
\putmorphism(\xpos,\ypos)(-1,-1)[\tokena``{#4}]{\height}{\arrowtypea}l%
\putvmorphism(\xpos,\ypos)[``{#5}]{\height}{\arrowtypeb}m%
\putmorphism(\xpos,\ypos)(1,-1)[``{#6}]{\height}{\arrowtypec}r%
}}
\def\putAtrianglepair{\@ifnextchar <{\putAtrianglepairp}{\putAtrianglepairp%
   <\arrowtypea`\arrowtypeb`\arrowtypec`\arrowtyped`\arrowtypee;\height>}}
\def\Atrianglepair{\@ifnextchar <{\Atrianglepairp}{\Atrianglepairp%
   <\arrowtypea`\arrowtypeb`\arrowtypec`\arrowtyped`\arrowtypee;\height>}}
\def\Atrianglepairp<#1>[#2;#3`#4`#5`#6`#7]{{
\settripairparms[#1]
\settokens[#2]
\width=\height                                
\diagram
\putAtrianglepairp                            
<\arrowtypea`\arrowtypeb`\arrowtypec`
\arrowtyped`\arrowtypee;\height>
(0,0)[{#2};#3`#4`#5`#6`{#7}]
\enddiagram
}}
\def\putVtrianglep<#1>(#2,#3)[#4`#5`#6;#7`#8`#9]{{%
\settriparms[#1]%
\xpos=#2 \ypos=#3 \advance\ypos by\height {\multiply\height by2
\puthmorphism(\xpos,\ypos)[#4`#5`{#7}]{\height}{\arrowtypea}a}%
\putmorphism(\xpos,\ypos)(1,-1)[`#6`{#8}]{\height}{\arrowtypeb}l%
\advance\xpos by\height \advance\xpos by\height
\putmorphism(\xpos,\ypos)(-1,-1)[``{#9}]{\height}{\arrowtypec}r%
}}
\def\putVtriangle{\@ifnextchar <{\putVtrianglep}{\putVtrianglep
   <\arrowtypea`\arrowtypeb`\arrowtypec;\height>}}
\def\Vtriangle{\@ifnextchar <{\Vtrianglep}{\Vtrianglep
   <\arrowtypea`\arrowtypeb`\arrowtypec;\height>}}
\def\Vtrianglep<#1>[#2`#3`#4;#5`#6`#7]{{
\settriparms[#1]
\width=\height                                 
\diagram
\putVtrianglep<\arrowtypea`\arrowtypeb`
\arrowtypec;\height>
(0,0)[#2`#3`#4;#5`#6`{#7}]
\enddiagram
}}
\def\putVtrianglepairp<#1>(#2)[#3;#4`#5`#6`#7`#8]{{
\settripairparms[#1]%
\setpos(#2)%
\settokens[#3]%
\advance\ypos by\height
\putmorphism(\xpos,\ypos)(1,-1)[`\tokend`{#6}]{\height}{\arrowtypec}l%
\puthmorphism(\xpos,\ypos)[\tokena`\tokenb`{#4}]{\height}{\arrowtypea}a%
\advance\xpos by\height
\puthmorphism(\xpos,\ypos)[\phantom{\tokenb}`\tokenc`{#5}]%
{\height}{\arrowtypeb}a%
\putvmorphism(\xpos,\ypos)[``{#7}]{\height}{\arrowtyped}m%
\advance\xpos by\height
\putmorphism(\xpos,\ypos)(-1,-1)[``{#8}]{\height}{\arrowtypee}r%
}}
\def\putVtrianglepair{\@ifnextchar <{\putVtrianglepairp}{\putVtrianglepairp%
    <\arrowtypea`\arrowtypeb`\arrowtypec`\arrowtyped`\arrowtypee;\height>}}
\def\Vtrianglepair{\@ifnextchar <{\Vtrianglepairp}{\Vtrianglepairp%
    <\arrowtypea`\arrowtypeb`\arrowtypec`\arrowtyped`\arrowtypee;\height>}}
\def\Vtrianglepairp<#1>[#2;#3`#4`#5`#6`#7]{{
\settripairparms[#1]
\settokens[#2]
\diagram
\putVtrianglepairp                             
<\arrowtypea`\arrowtypeb`\arrowtypec`
\arrowtyped`\arrowtypee;\height>
(0,0)[{#2};#3`#4`#5`#6`{#7}]
\enddiagram
}}

\def\putCtrianglep<#1>(#2,#3)[#4`#5`#6;#7`#8`#9]{{%
\settriparms[#1]%
\xpos=#2 \ypos=#3 \advance\ypos by\height
\putmorphism(\xpos,\ypos)(1,-1)[``{#9}]{\height}{\arrowtypec}l%
\advance\xpos by\height \advance\ypos by\height
\putmorphism(\xpos,\ypos)(-1,-1)[#4`#5`{#7}]{\height}{\arrowtypea}l%
{\multiply\height by 2
\putvmorphism(\xpos,\ypos)[`#6`{#8}]{\height}{\arrowtypeb}r}%
}}
\def\putCtriangle{\@ifnextchar <{\putCtrianglep}{\putCtrianglep
    <\arrowtypea`\arrowtypeb`\arrowtypec;\height>}}
\def\Ctriangle{\@ifnextchar <{\Ctrianglep}{\Ctrianglep
    <\arrowtypea`\arrowtypeb`\arrowtypec;\height>}}
\def\Ctrianglep<#1>[#2`#3`#4;#5`#6`#7]{{
\settriparms[#1]
\width=\height                               
\diagram
\putCtrianglep<\arrowtypea`\arrowtypeb`
\arrowtypec;\height>
(0,0)[#2`#3`#4;#5`#6`{#7}]
\enddiagram
}}                                           
\def\putDtrianglep<#1>(#2,#3)[#4`#5`#6;#7`#8`#9]{{%
\settriparms[#1]%
\xpos=#2 \ypos=#3 \advance\xpos by\height \advance\ypos by\height
\putmorphism(\xpos,\ypos)(-1,-1)[``{#9}]{\height}{\arrowtypec}r%
\advance\xpos by-\height \advance\ypos by\height
\putmorphism(\xpos,\ypos)(1,-1)[`#5`{#8}]{\height}{\arrowtypeb}r%
{\multiply\height by 2
\putvmorphism(\xpos,\ypos)[#4`#6`{#7}]{\height}{\arrowtypea}l}%
}}
\def\putDtriangle{\@ifnextchar <{\putDtrianglep}{\putDtrianglep
    <\arrowtypea`\arrowtypeb`\arrowtypec;\height>}}
\def\Dtriangle{\@ifnextchar <{\Dtrianglep}{\Dtrianglep
   <\arrowtypea`\arrowtypeb`\arrowtypec;\height>}}
\def\Dtrianglep<#1>[#2`#3`#4;#5`#6`#7]{{
\settriparms[#1]
\width=\height                              
\diagram
\putDtrianglep<\arrowtypea`\arrowtypeb`
\arrowtypec;\height>
(0,0)[#2`#3`#4;#5`#6`{#7}]
\enddiagram
}}                                          
\def\setrecparms[#1`#2]{\width=#1 \height=#2}%
\def\recursep<#1`#2>[#3;#4`#5`#6`#7`#8]{{%
\width=#1 \height=#2 \settokens[#3]
\settowidth{\tempdimen}{$\tokena$} \ifdim\tempdimen=0pt
  \savebox{\tempboxa}{\hbox{$\tokenb$}}%
  \savebox{\tempboxb}{\hbox{$\tokend$}}%
  \savebox{\tempboxc}{\hbox{$#6$}}%
\else
  \savebox{\tempboxa}{\hbox{$\hbox{$\tokena$}\times\hbox{$\tokenb$}$}}%
  \savebox{\tempboxb}{\hbox{$\hbox{$\tokena$}\times\hbox{$\tokend$}$}}%
  \savebox{\tempboxc}{\hbox{$\hbox{$\tokena$}\times\hbox{$#6$}$}}%
\fi \ypos=\height \divide\ypos by 2 \xpos=\ypos \advance\xpos by
\width \bfig
\putCtrianglep<-1`1`1;\ypos>(0,0)[`\tokenc`;#5`#6`{#7}]%
\puthmorphism(\ypos,0)[\tokend`\usebox{\tempboxb}`{#8}]{\width}{-1}b%
\puthmorphism(\ypos,\height)[\tokenb`\usebox{\tempboxa}`{#4}]{\width}{-1}a%
\advance\ypos by \width
\putvmorphism(\ypos,\height)[``\usebox{\tempboxc}]{\height}1r%
\efig }}
\def\recurse{\@ifnextchar <{\recursep}{\recursep<\width`\height>}}
\def\puttwohmorphisms(#1,#2)[#3`#4;#5`#6]#7#8#9{{%
\puthmorphism(#1,#2)[#3`#4`]{#7}0a \ypos=#2 \advance\ypos by 20
\puthmorphism(#1,\ypos)[\phantom{#3}`\phantom{#4}`#5]{#7}{#8}a
\advance\ypos by -40
\puthmorphism(#1,\ypos)[\phantom{#3}`\phantom{#4}`#6]{#7}{#9}b }}
\def\puttwovmorphisms(#1,#2)[#3`#4;#5`#6]#7#8#9{{%
\putvmorphism(#1,#2)[#3`#4`]{#7}0a \xpos=#1 \advance\xpos by -20
\putvmorphism(\xpos,#2)[\phantom{#3}`\phantom{#4}`#5]{#7}{#8}l
\advance\xpos by 40
\putvmorphism(\xpos,#2)[\phantom{#3}`\phantom{#4}`#6]{#7}{#9}r }}
\def\puthcoequalizer(#1)[#2`#3`#4;#5`#6`#7]#8#9{{%
\setpos(#1)%
\puttwohmorphisms(\xpos,\ypos)[#2`#3;#5`#6]{#8}11%
\advance\xpos by #8
\puthmorphism(\xpos,\ypos)[\phantom{#3}`#4`#7]{#8}1{#9} }}
\def\putvcoequalizer(#1)[#2`#3`#4;#5`#6`#7]#8#9{{%
\setpos(#1)%
\puttwovmorphisms(\xpos,\ypos)[#2`#3;#5`#6]{#8}11%
\advance\ypos by -#8
\putvmorphism(\xpos,\ypos)[\phantom{#3}`#4`#7]{#8}1{#9} }}
\def\putthreehmorphisms(#1)[#2`#3;#4`#5`#6]#7(#8)#9{{%
\setpos(#1) \settypes(#8)
\if a#9 %
     \vertsize{\tempcounta}{#5}%
     \vertsize{\tempcountb}{#6}%
     \ifnum \tempcounta<\tempcountb \tempcounta=\tempcountb \fi
\else
     \vertsize{\tempcounta}{#4}%
     \vertsize{\tempcountb}{#5}%
     \ifnum \tempcounta<\tempcountb \tempcounta=\tempcountb \fi
\fi \advance \tempcounta by 60
\puthmorphism(\xpos,\ypos)[#2`#3`#5]{#7}{\arrowtypeb}{#9}
\advance\ypos by \tempcounta
\puthmorphism(\xpos,\ypos)[\phantom{#2}`\phantom{#3}`#4]{#7}{\arrowtypea}{#9}
\advance\ypos by -\tempcounta \advance\ypos by -\tempcounta
\puthmorphism(\xpos,\ypos)[\phantom{#2}`\phantom{#3}`#6]{#7}{\arrowtypec}{#9}
}}
\def\setarrowtoks[#1`#2`#3`#4`#5`#6]{%
\def\toka{#1}
\def\tokb{#2}
\def\tokc{#3}
\def\tokd{#4}
\def\toke{#5}
\def\tokf{#6}
}
\def\hex{\@ifnextchar <{\hexp}{\hexp<1000`400>}}
\def\hexp<#1`#2>[#3`#4`#5`#6`#7`#8;#9]{%
\setarrowtoks[#9] \yext=#2 \advance \yext by #2 \xext=#1
\advance\xext by \yext \bfig
\putCtriangle<-1`0`1;#2>(0,0)[`#5`;\tokb``\tokd] \xext=#1
\yext=#2 \advance \yext by #2
\putsquare<1`0`0`1;\xext`\yext>(#2,0)[#3`#4`#7`#8;\toka```\tokf]
\advance \xext by #2
\putDtriangle<0`1`-1;#2>(\xext,0)[`#6`;`\tokc`\toke] \efig }
\begin{document}
\newtheorem{theorem}{Theorem}[section]
\newtheorem{lemma}[theorem]{Lemma}
\newtheorem{corollary}[theorem]{Corollary}
\newtheorem{conjecture}[theorem]{Conjecture}
\newtheorem{remark}[theorem]{Remark}
\newtheorem{condition}[theorem]{Condition}{\it}{\rm}
\newtheorem{definition}[theorem]{Definition}
\newtheorem{problem}[theorem]{Problem}
\newtheorem{example}[theorem]{Example}
\newtheorem{proposition}[theorem]{Proposition}
\newcommand{\cA}{{\mathcal A} }
\newcommand{\cB}{{\mathcal B} }
\newcommand{\cD}{{\mathcal D} }
\newcommand{\cE}{{\mathcal E} }
\newcommand{\cF}{{\mathcal F} }
\newcommand{\cG}{{\mathcal G} }
\newcommand{\cH}{{\mathcal H} }
\newcommand{\cI}{{\mathcal I} }
\newcommand{\cJ}{{\mathcal J} }
\newcommand{\cK}{{\mathcal K} }
\newcommand{\cL}{{\mathcal L} }
\newcommand{\cM}{{\mathcal M} }
\newcommand{\cN}{{\mathcal N} }
\newcommand{\cO}{{\mathcal O} }
\newcommand{\cP}{{\mathcal P} }
\newcommand{\cQ}{{\mathcal Q} }
\newcommand{\cS}{{\mathcal S} }
\newcommand{\cT}{{\mathcal T} }
\newcommand{\cW}{{\mathcal W} }
\newcommand{\cX}{{\mathcal X} }
\newcommand{\cY}{{\mathcal Y} }
\newcommand{\cZ}{{\mathcal Z} }
\newcommand{\bD}{{\bf D} }
\newcommand{\imp}{{\Rightarrow}}
\newcommand{\wt}{\widetilde}
\newcommand{\wh}{\widehat}
\newcommand{\Hom}{{\rm Hom}}
\def\ol#1{{\overline{#1}}}
\def\psh{{plurisubharmonic}}
\title{Dynamical construction of K\"{a}hler-Einstein metrics on  bounded pseudoconvex domains}
\author{Hajime Tsuji}
\maketitle
\begin{abstract}
\noindent We shall prove that the complete K\"{a}hler-Einstein metric on a bounded strongly pseudoconvex domain with $C^{\infty}$-boundary is 
the normalized  limit of a sequence of Bergman metrics. This is a noncompact version of \cite[p.110,Theorem 1,2]{tu8}. \end{abstract}
\section{Introduction}
Let $\Omega$ be a bounded strongly pseudoconvex domain in $\mathbb{C}^{n}$ with $C^{\infty}$-boundary $\partial\Omega$.  There has been constructed several canonical metrics on $\Omega$. 

Among them, the Bergman metric and the K\"{a}hler-Einstein metric are important. 
First the Bergman metric is constructed as follows.   Let $A^{2}(\Omega)$ be the space of $L^{2}$-holomorphic $n$-forms on $\Omega$, i.e., 
\[
A^{2}(\Omega):= \left\{\eta\in H^{0}(\Omega,\mathcal{O}_{\Omega}(K_{\Omega}))
\left|\int_{\Omega}|\eta|^{2} < \infty\right.\right\},
\]
where $|\eta|^{2}:= (\sqrt{-1})^{n^{2}}\eta\wedge\bar{\eta}$.
$A^{2}(\Omega)$ has a structure of a Hilbert space with repsect to the inner product:
\[
(\sigma,\tau):= (\sqrt{-1})^{n^{2}}\int_{\Omega}\sigma\wedge\bar{\tau}\,\,\,\,\,(\sigma,\tau \in A^{2}(\Omega)).
\] 
Let $\{\sigma_{j}\}_{j=1}^{\infty}$ be a complete orthonormal basis of $A^{2}(\Omega)$.  We set 
\[
K(\Omega)(z):= \sum_{j=1}^{\infty}|\sigma_{j}(z)|^{2} \,\,\,\,(z\in \Omega)
\]
and call it (the diagonal part of) the Bergman kernel or the Bergman volume form on $\Omega$. 
Then 
\[
\omega_{B}:= \sqrt{-1}\partial\bar{\partial}\log K(\Omega) 
\] 
is a complete K\"{a}hler form on $\Omega$ (\cite{kr}), 
$K(\Omega)$ has the following extremal property:
\[
K(\Omega)(z) = \sup\{|\sigma|^{2}(z)| \sigma\in A^{2}(\Omega),\parallel\sigma\parallel = 1\}. 
\]

Another important complete K\"{a}hler metric on $\Omega$ is the complete K\"{a}hler-Einstein metric. 
In \cite{c-y}, Cheng and Yau constructed a complete K\"{a}hler-Einstein metric on a bounded strongly pseudoconvex domain in $\mathbb{C}^{n}$ with $C^{\infty}$-boundary $\partial\Omega$. More precisely they proved the following theorem:
\begin{theorem}(\cite{c-y})\label{KE}
Let $\Omega$ be a bounded strongly pseudoconvex domain in $\mathbb{C}^{n}$ with $C^{\infty}$-boundary. 

Then there exists a complete $C^{\infty}$-K\"{a}hler-Einstein form $\omega_{E}$ on $\Omega$ such that 
\[
-\mbox{\em Ric}(\omega_{E}) = \omega_{E}
\]
holds on $\Omega$, where 
\[
\mbox{\em Ric}(\omega_{E}):= -\sqrt{-1}\partial\bar{\partial}\log\det(g_{i\bar{j}}).
\]
Here $(g_{i\bar{j}})$ is the hermitian matrix defind by 
\[
\omega_{E} = \frac{\sqrt{-1}}{2}\sum_{i,j=1}^{n}g_{i\bar{j}}\,dz_{i}\wedge dz_{j}.
\]\fbox{} 
\end{theorem}
We call the volume form 
\[
dV_{E}:= \frac{1}{n!}\,\omega_{E}^{n}
\]
associated with $\omega_{E}$ the K\"{a}hler-Einstein volume form on $\Omega$. \vspace{3mm}\\ 

The purpose of this paper is to relate the Bergman volume form $K(\Omega)$ and the K\"{a}hler-Einstein volume form $dV_{E}$ in terms of a dynamical system of 
Bergman kernels as in \cite{tu8}. \vspace{3mm}\\ 

First we set $K_{1}:= K(\Omega)$ and let $h_{1}:= K_{1}^{-1}$. Then $h_{1}$ is 
a $C^{\infty}$-hermitian metric on the canonical bundle $K_{\Omega}$ of $K_{\Omega}$ with strictly positive curvature. 

Suppose that we have already constructed $\{ K_{1},\cdots ,K_{m}\}$ and 
$\{ h_{1},\cdots ,h_{m}\}(h_{k}= K_{k}^{-1})$. 
Then we define 
\begin{equation}
K_{m+1}:= K(\Omega,(m+1)K_{\Omega},h_{m}),
\end{equation}
i.e.,
\[
K_{m+1} = \sum_{j=1}^{\infty}|\sigma_{j}^{(m+1)}|^{2}, 
\]
where $\{\sigma_{1}^{(m+1)},\cdots,\sigma_{k}^{(m+1)},\cdots\}$ is a complete orthonormal basis of the Hilbert space:
\[
A^{2}(\Omega,(m+1)K_{\Omega},h_{m})
:= \left\{\sigma\in H^{0}(\Omega,\mathcal{O}_{\Omega}((m+1)K_{\Omega}))\left|
\int_{\Omega}|\sigma|^{2}\cdot h_{m} < +\infty \right.\right\},
\]
with respect to the inner product:
\[
(\sigma,\tau)_{(m+1)}= (\sqrt{-1})^{n^{2}}\int_{\Omega}\sigma\wedge\bar{\tau}\cdot h_{m}. 
\]
And we set 
\[
h_{m+1}:= K_{m+1}^{-1}.
\]
Since $A^{2}(\Omega,(m+1)K_{\Omega},h_{m})$ is very ample by H\"{o}rmander's 
$L^{2}$-estimate of $\bar{\partial}$-operaters, we see that 
the dynamical system $\{ K_{m}\}_{m=1}^{\infty}$ is well defined. The dynamical system is essentially the same as in \cite{tu8}. The following is the main theorem in this paper. 

\begin{theorem}\label{main}
Let $\Omega$ be a bounded strongly pseudoconvex domain in $\mathbb{C}^{n}$ with $C^{\infty}$-boundary.  Let $\{ K_{m}\}_{m=1}^{\infty}$ be the dynamical system of Bergman kernels constructed as above. 

Then $\lim_{m\rightarrow\infty}\sqrt[m]{(m!)^{-n}K_{m}}$ exists in compact 
uniform topology and 
\[
\lim_{m\rightarrow\infty}\sqrt[m]{(m!)^{-n}K_{m}} = (2\pi)^{-n}dV_{E}
\]
holds. \fbox{}
\end{theorem}
The proof of Theorem \ref{main} is very similar to the proof of \cite[p.110,Theorem 1.2]{tu8}.  We just need to take care of the unifomity of the estimates, since $\Omega$ is noncompact. 

\section{Proof of Theorem \ref{main}}

To ensure the uniformity of the estimate, we shall introduce the following notion. 

\begin{definition}\label{bounded}  Let $M$ be a complete K\"{a}hler manifold.  We say that $M$ has bounded geometry of order $\ell$, if $M$ admits a covering of holomorphic 
coodinate charts $\{(V,(v^{1},\cdots ,v^{n}))\}$ and positive numbers $R,c,\mathcal{A}_{1},\cdots ,\mathcal{A}_{\ell}$ such that 
\begin{enumerate}
\item For any $x_{0}\in M$, there is a coordinate chart $(V,(v^{1},\cdots,v^{n}))$ with $x_{0}\in V$ and with repsect to the Euclidean distance $d$ defined by $v^{i}$-coordinates, $d(x_{0},\partial V) > \sqrt{n}\cdot R$ holds. 
\item If $(g_{i\bar{j}})$ denotes the metric tensor with repsect to 
$(V,(v^{1},\cdots ,v^{n}))$, then $(g_{i\bar{j}})$ is $C^{\ell}$ and 
\[
\frac{1}{C}\,(\delta_{ij}) < (g_{i\bar{j}}) < C(\delta_{ij})
\]
holds and for any multi index $\alpha,\beta$ with $|\alpha| + |\beta| \leqq \ell$, we have 
\[
\left|\frac{\partial^{|\alpha|+|\beta|}}{\partial v^{\alpha}\partial\bar{v}^{\beta}}\,g_{i\bar{j}}\right| \leqq \mathcal{A}_{|\alpha|+|\beta|}
\]
holds. \fbox{}
\end{enumerate}
\end{definition}

\begin{theorem}(\cite{c-y})\label{BG}
Let $\Omega$ be a bounded strongly pseudoconvex domain in $\mathbb{C}^{n}$ with $C^{\infty}$-boundary.  Let $\omega_{E}$ be the complete K\"{a}hler-Einstein form $\omega_{E}$ such that $-\mbox{\em Ric}(\omega_{E}) = \omega_{E}$ on $\Omega$ as in Theorem \ref{KE}.  

Then $(\Omega,\omega_{E})$ has bounded geometry of $\infty$-order. \fbox{} 
\end{theorem}

Let $(\Omega,\omega_{E})$ be a bounded strongly pseudoconvex domain with $C^{\infty}$-boundary. Let $\{ K_{m}\}_{m=1}^{\infty}$ be the dynamical system of Bergman kernels as in Section 1.  Hereafter we shall estimate $\{ K_{m}\}_{m=1}^{\infty}$. 

\subsection{Upper estimate}

By Theorem \ref{BG}, $(\Omega,\omega_{E})$ has bounded geometry of $\infty$-order. Hence there exist positive numbers $R,c,\mathcal{A}_{1},\cdots,\mathcal{A}_{\ell},\cdots$ and a covering of holomorphic coordinate charts $\{(V,(v_{1},\cdots,v_{n}))$ satisfying the conditions in Definition \ref{bounded}. Let $x_{0}\in \Omega$ be an arbitrary point. Then since $(\Omega,\omega_{E})$ has a bounded geometry of $\infty$-order, there exists a local coordinate $(U,(z_{1},\cdots,z_{n}))$ such that 
\begin{enumerate}
\item[(1)] $U$ is biholomorphic to the unit open polydisc $\Delta^{n}(R)$ with 
center $O$ via $(z_{1},\cdots,z_{n})$,
\item[(2)] $g_{i\bar{j}} = \delta_{ij} + O(\parallel z\parallel^{2})$,
\item[(3)] 
\begin{equation}\label{det}
\det(g_{i\bar{j}}) = \left(\prod_{i=1}^{n}\left(1-\frac{1}{2}|z_{i}|^{2}\right)\right)^{-1}
+ O(\parallel z\parallel^{3}). 
\end{equation}
\end{enumerate}
First we note that 
\[
K_{1}(x_{0}) \leqq K(\Delta^{n}(R))(O) = \frac{1}{2^{n}(\pi R^{2})^{n}}|dz_{1}\wedge\cdots \wedge dz_{n}|^{2}
\]
holds. Hence there exists a positive constant $C_{1,+}$ such that  
\begin{equation}
K_{1} \leqq C_{1,+}\cdot dV_{E}
\end{equation}
holds on $\Omega$. 

Now  we proceed by induction on $m$.  Suppose that for some $m\geqq 1$, there exists a positive constant $C_{m,+}$ such that 
\begin{equation}
K_{m} \leqq C_{m,+}\cdot (dV_{E})^{m}
\end{equation}
holds on $\Omega$. We note that by the extremal property of the Bergman kernel
\begin{equation}
K_{m+1}(x) = \sup\left\{|\sigma|^{2}(x)\left|\sigma\in \Gamma(\Omega,\mathcal{O}_{\Omega}((m+1)K_{\Omega})),\int_{\Omega}|\sigma|^{2}\cdot h_{m} = 1\right.\right\}
\end{equation}
holds. Then by the induction hypothesis, we have that for every $r \leqq R$, 
\begin{equation}
K_{m+1}(x_{0})\leqq C_{m}\cdot K(B(x_{0},r),(m+1)K_{\Omega},dV_{E}^{-m})
\end{equation}
holds. By the Taylor expansion (\ref{det}) of $\det(g_{i\bar{j}})$, we have that there exists a positive constant $c < 1$ independent of $m$ and a positive function $\delta(r)$ of $r$ such that 
\begin{eqnarray}
K(B(x_{0},r),(m+1)K_{\Omega},dV_{E}^{-m})(x_{0}) \leqq  \hspace{30mm} \\   
\left(\frac{m+1}{2\pi}\right)^{n}\!\!\!\cdot(1-c^{m+1})^{-1}(1+\delta(r))\cdot (2^{-n}|dz_{1}\wedge\cdots\wedge dz_{n}|^{2})^{m+1} \nonumber
\end{eqnarray}
and $\lim_{r\downarrow 0}\delta(r) = 0$ hold.   
Here $c$ corresponds to the fact that for a positive number $\rho < 1$ 
\[
\frac{\sqrt{-1}}{2}\int_{\Delta(\rho)}\left(1-\frac{1}{2}|t|^{2}\right)^{m}dt\wedge d\bar{t} = \frac{2\pi}{m+1}\left(1-\left(1-\frac{1}{2}\rho^{2}\right)^{m+1}\right)
\]
holds, where $\Delta(\rho) = \{ t\in \mathbb{C}|\,\, |t| < \rho\}$.  And $\delta(r)$ corresponds to the Taylor expansion $(\ref{det})$. 
This implies that 
\begin{equation}
\limsup_{m\to\infty}\sqrt[m]{(m!)^{-n}K_{m}} \leqq (1 + \delta(r))(2\pi)^{-n}\cdot 2^{-n}|dz_{1}\wedge\cdots\wedge dz_{n}|^{2}
\end{equation} 
holds. 
Since the estimate is independent of $r$, letting $r$ tend to $0$, noting 
$\lim_{r\downarrow 0}\delta(r) = 0$, we have the following lemma. 

\begin{lemma}\label{upper}
\[
\limsup_{m\to\infty}\sqrt[m]{(m!)^{-n}K_{m}} \leqq (2\pi)^{-n}dV_{E}
\]
holds on $\Omega$. \fbox{}
\end{lemma}

\subsection{Lower estimate}

Now we shall estimate $\{ K_{m}\}_{m=1}^{\infty}$ from below. 
The lower estimate is also similar to the one in \cite[Section 3.2]{tu8}.  The only difference is the use of the fact that $(\Omega,\omega_{E})$ has a bounded geometry of 
$\infty$-order.

Take a defining function $\rho$ of $\Omega$ and set $\rho_{\sharp}(z_{0}, z) = |z_{0}|^{2}\rho(z)$, which is a defining function of $\mathbb{C}^{*}\times \Omega$. The strictly
pseudoconvexity of $\Omega$ then ensures that
\[
g[\rho] =\sum^{n}_{j,k=0}
\frac{\partial^{2}\rho_{\sharp}}{\partial z_{j}\partial\bar{z}_{k}}dz_{j}d\bar{z}_{k}
\]
is a Lorentz-K\"{a}hler metric in a neighborhood of $\mathbb{C}^{*}\times\partial\Omega$ in $\mathbb{C}^{*}\times \overline{\Omega}$. 
We shall recall the following theorem.

\begin{theorem}\label{BA} The Bergman kenel of a strongly pseudoconvex domain $\Omega$ with $C^{\infty}$-boundary has the following singularity at the boundary $\partial\Omega$. 

Let $r$ be a defining function of $\Omega$ satisfying 
\[
J[r] = 1 + O^{n+1}(\partial\Omega),
\]
where
\[
J[r] = (-1)^{n}\det\left(\begin{array}{ll} r & \partial r/\partial z_{j} \\
\partial r/\partial\bar{z}_{k} & \partial^{2}r/\partial z_{j}\partial\bar{z}_{k}\end{array}\right) (j,k = 1,\cdots,n)
\]
and $f =  O^{n+1}(\partial\Omega)$ stands for a function such that $f/r^{n}\in C^{\infty}(\overline{\Omega})$.  Then 
\[
K(\Omega) = r^{-n-1}(c_{n} + c_{n}^{\prime}\parallel R\parallel^{2}r^{2} + \cdots),
\]
wnere $c_{n} = n!/\pi^{n}$ , $c_{n}^{\prime} = (n-2)!/(24\pi^{n})$ and 
$R$ is the curvature of the Lorenz K\"{a}hler metric near the boundary. \fbox{} 
\end{theorem}
On the other hand, the complete K\"{a}hler-Einstein metric of \cite{c-y} is equivalent to the model metric constructed as follows.

Let 
$\Omega = \{\varphi < 0\}$.   We set 
\begin{equation}
\omega := \sqrt{-1}\partial\bar{\partial}(-\log (-\varphi)) 
= \frac{\sqrt{-1}}{2}\sum g_{i\bar{j}}dz_{i}\wedge d\bar{z}_{j}
\end{equation}
Then by the direct calculation, we have 
\begin{equation}\label{modelboundary}
\det (g_{i\bar{j}}) = \left(-\frac{1}{\varphi}\right)^{n+1}\!\!\!\!\!\det(\varphi_{i\bar{j}})(-\varphi + |d\varphi|^{2})
\end{equation}
holds (cf. \cite[p.510]{c-y}).  By the construction of $\omega_{E}$. we have that there exists a positive constant $C$ such that 
\begin{equation}\label{qi}
\frac{1}{C}\,\omega \leqq \omega_{E} \leqq C\omega
\end{equation}
holds on $\Omega$.  
Hence by Theorem \ref{BA}, (\ref{modelboundary}) and (\ref{qi}), we see that 
there exists a positive constant $C_{1,-}~{\prime}$ such that 
\begin{equation}\label{m=1}
K_{1} \geqq C_{1,-}^{\prime}\cdot dV_{E}
\end{equation}
holds on $\Omega$.  

Now we shall use the fact that $(\Omega,\omega_{E})$ has bounded geometry.  
Let $x_{0}$ be an arbitrary point and let $(U,(z_{1},\cdots,z_{n})$ be a holomorphic local coordinates satisfying following conditions:  
\begin{enumerate}
\item[(1)] $U$ is biholomorphic to the unit open polydisc $\Delta^{n}(R)$ with 
center $O$ via $(z_{1},\cdots,z_{n})$,
\item[(2)] $g_{i\bar{j}} = \delta_{ij} + O(\parallel z\parallel^{2})$,
\item[(3)] 
\[
\det(g_{i\bar{j}}) = \left(\prod_{i=1}^{n}\left(1- \frac{1}{2}|z_{i}|^{2}\right)\right)^{-1}
+ O(\parallel z\parallel^{3}). 
\]
Let $\rho$ be a $C^{\infty}$-function on $\Omega$ such that 
\begin{enumerate}
\item[(1)] $\rho\equiv 1$ on $B(x_{0},R/3)$,
\item[(2)] $0\leqq \rho \leqq 1$, 
\item[(3)] $\mbox{supp}\,\rho \subset\subset B(x_{0},R)$,
\item[(4)] $\rho$ is a function of the Euclidean distance from $x_{0}$ 
on $B(x_{0},R)$ with respect to $(z_{1},\cdots,z_{n})$, 
\item[(5)] $| d\rho| \leqq 3/R$. 
\end{enumerate}
\end{enumerate} \label{phi}
Let $\phi$ be the function defined by 
\begin{equation}
\phi := \alpha\cdot\log\frac{dV_{E}}{dV_{\mathbb{C}^{n}}}, 
\end{equation}
where $dV_{\mathbb{C}^{n}}$ denotes the Euclidean volume form on $\mathbb{C}^{n}$ and $\alpha > 0$ is a positive constant which will be specified later. 

Since $(\Omega,\omega_{E})$ has a bounded geometry, we see that if we take $\alpha$ sufficiently large, we may assume that for every $x_{0}\in \Omega$,
\begin{equation}
\phi_{x_{0}}:= \phi + n\rho\log |z|^{2}
\end{equation} 
is plurisuharmonic on $\Omega$.

Now we proceed by induction on $m$.  Suppose that there exists a positive constant $C_{m,-}$ such that 
\begin{equation}\label{ind}
K_{m} \geqq C_{m,-}\cdot (dV_{E})^{m}\cdot e^{-\phi} 
\end{equation} 
holds on $\Omega$.  For $m = 1$, this inequality certainly holds by (\ref{m=1}), if we take $C_{1,-}$ sufficiently small.  Suppose that for some $m\geqq 1$, (\ref{ind}) holds for 
some positive constant  $C_{m,-}$. 

Let $\sigma_{0}$ be a local holomorphic section of $(m+1)K_{\Omega}$ on $U$ defined by
\[
\sigma_{0}:= (dz_{1}\wedge\cdots\wedge dz_{n})^{\otimes m+1}. 
\]
We shall solve the $\bar{\partial}$-equation 
\begin{equation}
\bar{\partial}u = \bar{\partial}(\rho\cdot\sigma_{0}).  
\end{equation}
Then by the $L^{2}$-estimate for $\bar{\partial}$-operators, there exists a 
$C^{\infty}$-solution $u$ such that 
\begin{equation}
\int_{\Omega}|u|^{2}\cdot e^{-\phi_{x_{0}}}\cdot (dV_{E})^{-m}
\leqq \frac{1}{m}\int_{\Omega}|\bar{\partial}(\rho\cdot\sigma_{0})|^{2}\cdot e^{-\phi_{x_{0}}}(dV_{E})^{-m}.
\end{equation}  
Then 
\[
\sigma := \rho\cdot\sigma_{0} - u
\]
is a global holomorphic section of $(m+1)K_{\Omega}$ such that 
\begin{equation}
\left(\int_{\Omega}|\sigma|^{2}\cdot e^{-\phi}\cdot (dV_{E})^{-m}\right)^{\frac{1}{2}}
\leqq 
\left(1 - \frac{1}{\sqrt{m}}\right)\left(\int_{\Omega}|\rho\cdot\sigma_{0}|^{2}
\cdot(dV_{E})^{-m}\right)^{\frac{1}{2}}. 
\end{equation}
By the Taylor expansion 
\[
\det(g_{i\bar{j}}) = \left(\prod_{i=1}^{n}\left(1- \frac{1}{2}|z_{i}|^{2}\right)^{-1}\right)
+ O(\parallel z\parallel^{3}),  
\]
we see that 
\[
\int_{\Omega}|\rho\cdot\sigma_{0}|^{2}(dV_{E})^{-m}
\sim \left(\frac{2\pi}{m +1}\right)^{n}
\]
holds, in the sense that the ratio of the both sides tend to $1$ as $m$ goes to infinity. 
This implies that if we take $\alpha$ sufficiently large  
\begin{equation}
K_{m+1}(x_{0})\geqq  C_{m,-}\left(1 -\frac{1}{\sqrt{m}}\right)\cdot (2\pi)^{-n}(m+1)^{n}\cdot (dV_{E})^{m+1}(x_{0})\cdot e^{-\phi(x_{0})}
\end{equation}
holds.  Hence we obtain the following lemma. 
\begin{lemma}\label{lower}
\[
\liminf_{m\to\infty}\sqrt[m]{(m!)^{-n}K_{m}} \geqq (2\pi)^{-n}dV_{E}
\]
holds on $\Omega$. \fbox{}
\end{lemma}

Combining Lemmas \ref{upper} and \ref{lower}, we complete the proof of 
Theorem \ref{main}. 

\section{Applications}

Theorem \ref{KE} has been generalized as follows. 

\begin{theorem}(\cite{m-y})\label{M-Y}
Let $\Omega$ be a bounded pseudoconvex domain in $\mathbb{C}^{n}$. Then there exixts a unique complete $C^{\infty}$-K\"{a}hler form $\omega_{E}$ on $\Omega$ such that $-\mbox{\em Ric}(\omega_{E}) = \omega_{E}$ holds on $\Omega$. \fbox{}
\end{theorem}  

The construction of $\omega_{E}$ is as follows. By the assumption, there exists a $C^{\infty}$-strictly plurisubharmonic exhaustion function $\varphi :\Omega\to \mathbb{R}$.  

For $c\in \mathbb{R}$, we set $\Omega_{c}:= \{\varphi < c\}$. Then there exists a set $E$ of measure $0$ in $\mathbb{R}$ such that for every $c\in \mathbb{R}\backslash E$, $\Omega_{c}$ is a bounded strongly pseudoconvex domain in $\mathbb{C}^{n}$ with $C^{\infty}$-boundary.

Then by Theorem \ref{KE}, we have the canonical complete K\"{a}hler-Einstein form $\omega_{E,c}$ on $\Omega_{c}$ such that $-\mbox{Ric}(\omega_{E,c}) = \omega_{E,c}$ holds on $\Omega_{c}$. By Yau's Schwarz lemma (\cite{y2}), for 
$c < c^{\prime}, c,c^{\prime}\in \mathbb{R}\backslash E$, 
\[
\omega_{E,c^{\prime}}^{n} \leqq \omega_{E,c}^{n}
\]
holds on $\Omega_{c}$. 
Hence $\{\omega_{E,c}| c\in \mathbb{R}\backslash E\}$ is a monotone decreasing 
sequence in $c \in \mathbb{R}\backslash E$.  Hence the limit
\[
dV_{E} := \frac{1}{n!}\lim_{c\to\infty}\omega_{E,c}^{n}
\]
exists on $\Omega$. It is easy to see that $dV_{E}$ is $C^{\infty}$ and 
$\omega_{E}:= - \mbox{Ric}\,dV_{E}$ satisfies the equality:
$-\mbox{Ric}(\omega_{E}) = \omega_{E}$.   

This  construction has already been considered in \cite{c-y}.  But the completeness of $\omega_{E}$ was first proved in \cite{m-y}. 

\begin{theorem}\label{psh}
Let $\Omega$ be a bounded pseudoconvex domain in $\mathbb{C}^{n}\times \Delta$.
Let $dV_{s}(s\in \Delta)$ denote the K\"{a}hler-Einstein volume form on 
$\Omega_{s}:= \Omega \cap (\mathbb{C}^{n}\times\{s\})$.

Let $dV_{\Omega/\Delta}$ be the relative volume form on $\Omega$ defined by
\[
dV_{\Omega/\Delta}|_{\Omega_{s}} = dV_{s}. 
\]
Then 
\[
\sqrt{-1}\partial\bar{\partial}\log dV_{\Omega/\Delta} \geqq 0
\]
holds on $\Omega$, i.e., $dV_{\Omega/\Delta}^{-1}$ is a hermitian metric on 
$K_{\Omega/\Delta}$ with semipositive curvature.  \fbox{}
\end{theorem}
\begin{remark}
By the construction of the complete K\"{a}hler-Einstein metric in \cite{c-y} 
and the implicit function theorem, we see that $dV_{\Omega/\Delta}$ is 
a $C^{\infty}$-relative volume form. 
\end{remark}
{\em Proof of Theorem \ref{psh}}. First we shall assume that $\Omega$ is a strongly pseudoconvex domain with $C^{\infty}$-boundary. Let $\{ K_{m,s}\}_{s\in\Delta}$ be the dynamical system of Bergman kernels on $\Omega_{s}$ constructed as in the previous section. 
Let $h_{m}$ be the hermitian metric on $mK_{X/S}$ defined by 
\[
h_{m}|_{\Omega_{s}}= K_{m,s}^{-1}. 
\]
Then by Berndtsson's theorem (\cite{b1,be}), we see that 
\[
\sqrt{-1}\Theta(h_{m}) > 0
\]
holds on $\Omega$ by induction on $m$. 
Then by Theorem \ref{main}, we see that $dV_{\Omega/\Delta}^{-1}$ has 
semipositive curvature. 

If $\Omega$ is a general bounded pseudoconvex domain in $\mathbb{C}^{n}\times\Delta$, there exists a $C^{\infty}$-strictly plurisubharmonic exhaustion function $\varphi$ on $\Omega$.  Then for every $c\in \mathbb{R}$, $\Omega_{c}:= \{ \varphi < c\}$ is a bounded strongly pseudoconvex domain in $\mathbb{C}^{n}\times \Delta$.  Then the metric $dV_{\Omega_{c}/\Delta}^{-1}$ on $K_{\Omega_{c}/\Delta}$  has semipositive curvature. 
Let $dV_{c,s}$ denote the K\"{a}hler-Einstein volume form on 
$\Omega_{c,s}$.  Then $dV_{c,s}$ is monotone decreasing in $c\in \mathbb{R}$. 
Since the limit of plurisubharmonic function is also plurisubharmonic, we complete the proof of Theorem \ref{psh}. \fbox{}

\small{

}

\noindent Author's address\\
Hajime Tsuji\\
Department of Mathematics\\
Sophia University\\
7-1 Kioicho, Chiyoda-ku 102-8554\\
Japan \\
h-tsuji@h03.itscom.net 


\begin{thebibliography}{99}
\bibitem[B-E-G]{b-e-g}T. N. Bailey, M. G. Eastwood and C. R. Graham, Invariant theory for conformal and CR
geometry, Ann. of Math. {\bf 139} (1994), 491-552
\bibitem[Ber1]{b1}B. Berndtsson: Subharmonicity properties of Bergman kernels and some other functions associated with pseudoconvex domains, Ann. Institute Fourier,{\bf 56} no.6, (2006),1633-1662.
\bibitem[Ber2]{be}B. Berndtsson: Curvature of vector bundles associated to holomorphic fibrations, Ann. of Math. {\bf 169} (2009), 531-560.
\bibitem[B-P]{b-p}B. Berndtsson and M. Paun : 
Bergman kernels and the pseudoeffectivity of relative canonical bundles, 
 Duke Math. J. {\bf 145}(2008), 341-378. .
\bibitem[C-Y]{c-y} S.-Y. Cheng and S.-T. Yau: On the Existence of a Complete K\"{a}hler Metric on Non-Compact Complex Manifolds and the Regularity of Fefferman's Equation, Communucation of Pure and Applied Mathematics, Ser.A {\bf 33}(1982), 507-544.   
\bibitem[F]{f} C. Fefferman, Parabolic invariant theory in complex analysis, Adv. in Math. {\bf 31} (1979), 131-262

\bibitem[Kr]{kr} S. Krantz : Function theory of several complex variables, 
John Wiley and Sons (1982).  
\bibitem[L]{l} P. Lelong : Fonctions Plurisousharmoniques et Formes 
Differentielles Positives, Gordon and Breach (1968).
\bibitem[M-Y]{m-y} N. Mok and S.-T. Yau : Completeness of the K\"{a}hler-Einstein metric on bounded
domains and the characterization of domains of holomorphy by curvature conditions, The Mathematical Heritage of Henri Poincar\'{e}, Proc. of Symp. in Pure Math., vol. {\bf 39}, Amer. Math. Soc., (1983), pp. 41-59.
\bibitem[T]{tu8} H. Tsuji: Dynamical construction of K\"{a}hler-Einstein metrics,  Nagoya Math. J. {\bf 199} (2010), 107-122.
\bibitem[Y]{y2} S.-T. Yau: A general Schwarz lemma for K\"{a}hler manifolds, Amer. J. of Math. {\bf 100} (1978), 197-203.

\end{thebibliography}
\end{document}